\documentclass[12pt]{article}

\usepackage[utf8]{inputenc}
\usepackage[english]{babel}
\usepackage{amssymb,amsmath}
\usepackage{amsthm}
\usepackage{float}
\usepackage{setspace,lipsum}
\usepackage{amsfonts}
\usepackage{graphicx,caption}
\usepackage[margin=2.4cm]{geometry}
\usepackage{enumerate}
\usepackage{enumitem}
\usepackage{color}
\usepackage[noadjust]{cite}
\usepackage{tikz}
\usepackage{mathtools}
\usepackage{multirow}
\usepackage{array}
\usepackage[section]{placeins}
\usepackage{xfrac}
\usepackage{xspace}
\usepackage{ifthen}
\usepackage{longtable}
\usepackage{subcaption}

\newtheorem{thm}{Theorem}
\newtheorem{conj}[thm]{Conjecture}

\theoremstyle{definition}

\title{PHOEG: an online tool for discovery and education in extremal graph theory}

\author{
    S\'ebastien Bonte,
	Gauvain Devillez,
	Valentin Dusollier,
	Hadrien M\'elot\textsuperscript{1}\\[3mm]
	\footnotesize  Computer Science Department - Algorithms Lab\\[-2mm]
	\footnotesize University of Mons, Mons, Belgium\\[-2mm]
	\footnotesize \textsuperscript{1}Corresponding author: hadrien.melot@umons.ac.be\\[3mm]
	\footnotesize
}

\date{\today}

\usepackage{hyperref}

\begin{document}

\maketitle

\vspace*{0.3cm}

\small
\noindent
\textbf{Abstract.}

\emph{
Extremal Graph Theory heavily relies on exploring bounds and inequalities between graph invariants, a task complicated by the rapid combinatorial explosion of graphs. Various tools have been developed to assist researchers in navigating this complexity, yet they typically rely on heuristic, probabilistic, or non-exhaustive methods, trading exactness for scalability. PHOEG takes a different stance: rather than approximating, it commits to an exact approach.}

\emph{PHOEG is an interactive online tool (\url{https://phoeg.umons.ac.be}) designed to assist researchers and educators in graph theory. Building upon the exact geometrical approach of its predecessor, GraPHedron, PHOEG embeds graphs into a two-dimensional invariant space and computes their convex hull, where facets represent inequalities and vertices correspond to extremal graphs. PHOEG modernizes and expands this approach by offering a comprehensive web interface and API, backed by an extensive database of pairwise non-isomorphic graphs including all graphs up to order 10. Users can intuitively define invariant spaces by selecting a pair of invariants, apply constraints and colorations, visualize resulting convex polytopes, and seamlessly inspect the corresponding drawn graphs.}

\emph{In this paper, we detail the software architecture and new web-based features of PHOEG. Furthermore, we demonstrate its practical value in two primary contexts: in research, by illustrating its ability to quickly identify conjectures or counterexamples to conjectures, and in education, by detailing its integration into university-level coursework to foster student discovery of classical graph theory principles. Finally, this paper serves as a brief survey of the extremal results and conjectures established over the past two decades using this geometric approach.
}

\vspace*{0.2cm}
\noindent
\textbf{Keywords:} chemical graphs, degree-based topological index, extremal graphs.

\vspace*{0.2cm}
\hrule

\normalsize

\section{Introduction}

Extremal Graph Theory studies bounds on graph invariants, which are values preserved under isomorphism, such as order, size, chromatic number, or average distance. These problems are typically subject to constraints like fixing the value of some invariants. Finding these bounds, and the extremal graphs that realize them, is challenging: the sheer number of graphs grows extremely rapidly, exceeding a hundred billion for graphs of order 12 alone.

To tackle this complexity, several computational tools have been developed over the years. Early efforts include Graph by Cvetkovic et al.~\cite{cvetkovic1981discussing} and its successor newGRAPH by Brankov et al.~\cite{stevanovic2003invitation}. Graffiti~\cite{fajtlowicz1988conjectures}, developed by Fajtlowicz in 1988, used heuristics and pre-computed data to generate thousands of conjectures automatically. AutoGraphiX~\cite{caporossi2000variable}, by Caporossi and Hansen in 2000, applied variable neighborhood search to identify extremal graph candidates.  Concurrently, databases like House of Graphs~\cite{brinkmann2013house, coolsaet2023house, devillez2016house} emerged, offering a searchable repository of interesting graphs allowing researchers to query specific invariant values or search for counterexamples. These tools have proven valuable, but they share a common characteristic: they are not exhaustive. Whether through heuristics, sampling, or selective search, they sacrifice exactness for scalability.

GraPHedron~\cite{melot2008facet}, introduced by Mélot in 2008, took a deliberately different stance: rather than approximating, it committed to an exact approach. It embedded all graphs up to a given order into an invariant space and computed their convex hull, interpreting facets as inequalities between invariants and vertices as extremal graphs. Historically, GraPHedron was instrumental in the early days of House of Graphs, as it was used to compute and introduce the very first graphs into the database. However, unlike some of its contemporaries, GraPHedron itself was never updated — motivating the development of its successor, PHOEG (a recursive acronym for PHOEG Helps to Obtain Extremal Graphs).

First introduced in 2019~\cite{devillez2019phoeg}, PHOEG embraces this same goal of exactness. The trade-off is explicit and assumed: PHOEG operates on all pairwise non-isomorphic graphs up to order 10, the practical frontier imposed by combinatorial growth. It combines a database of graphs enabling fast queries and computations with a proof-assistance module called TransProof. Since its initial release, we have significantly expanded PHOEG  to democratize its access (\url{https://phoeg.umons.ac.be}). Through extensive use in both research and educational settings, we refined the tool and developed a comprehensive web interface and API. 

This paper presents these modernizations and illustrates them through a detailed walkthrough of the platform. Section~\ref{sec:geometrical_approach} briefly recalls the underlying geometrical approach PHOEG is built on. Section~\ref{sec:interface_api} details the new web interface functionalities and the technologies used to develop them. Section~\ref{sec:phoeg_research} serves as a mini-survey of the extremal results obtained over the years using this geometric approach, and illustrates PHOEG's practical applications in research through concrete examples. Section~\ref{sec:phoeg_education} details its use in educational settings. Finally, in Section~\ref{sec:conclusion} we conclude and discuss future directions.

\section{Geometrical approach}
\label{sec:geometrical_approach}

This section explains the geometrical approach underlying PHOEG and introduced in GraPHedron~\cite{melot2008facet}. We first introduce a few basic mathematical definitions necessary to understand this methodology. Throughout this paper, we consider simple, undirected graphs $G = (V, E)$, where $V$ is the set of vertices and $E$ is the set of edges. The \emph{order} of a graph $G$ is the number of its vertices, denoted by $n = |V|$, and its \emph{size} is the number of its edges, denoted by $m = |E|$. The \emph{degree} of a vertex $v \in V$, denoted by $d(v)$, is the number of edges incident to $v$. For standard notations of graph theory that are not defined here, we refer the reader to Diestel~\cite{Diestel}.

As stated by Hansen et al.~\cite{hansen2005forms}, many results in extremal graph theory are expressed in terms of inequalities between graph invariants. This observation guided the design of GraPHedron and remains the core principle in PHOEG.

The geometrical approach used in PHOEG places simple, undirected graphs into a two-dimensional space where each dimension corresponds to a specific graph invariant. Each graph of a given order is represented as a point whose coordinates are the values of the two chosen invariants evaluated on that graph. It is important to note that a single point in this space may represent several non-isomorphic graphs, as they can share the exact same values for the chosen pair of invariants. The convex hull of all these points forms a convex polytope.  Recall that the \emph{convex hull} of a finite set of points in the plane is the smallest convex set enclosing them, and that its boundary consists of \emph{facets}, the edges of the polygon in two dimensions. Consequently, each facet of this convex polytope corresponds to a linear inequality between the invariants (typically of the form $aX + bY \leq c$), and the vertices represent the extremal graphs that tightly achieve these bounds.

\begin{figure}[!ht]
\begin{center}
\includegraphics[width=0.9\textwidth]{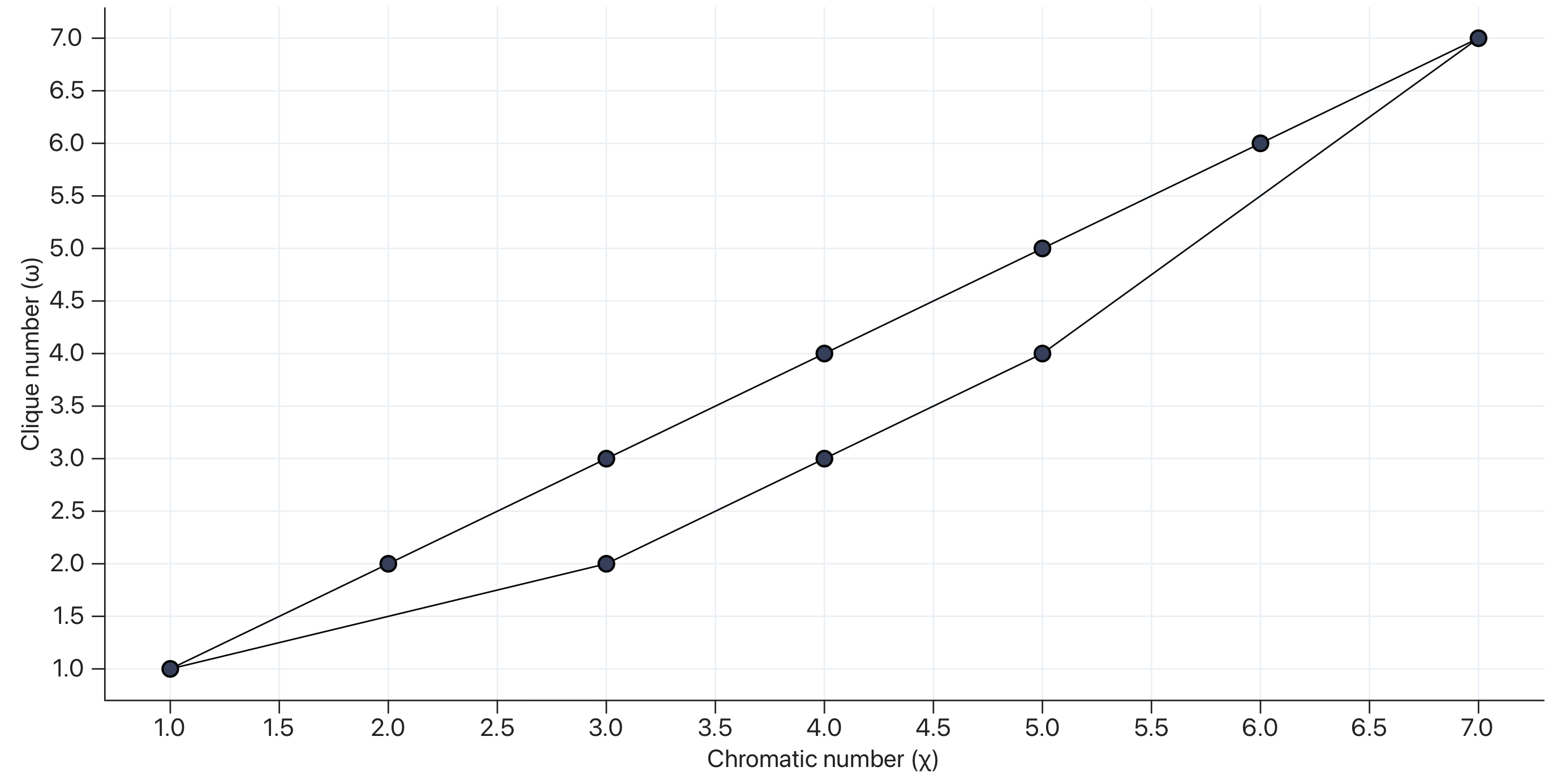}
\caption{Convex polytope formed by all graphs of order 7, with the chromatic number and the clique number used as dimensions}\label{fig:geom_example}
\end{center}
\end{figure}

An example of such a polytope is shown in Figure~\ref{fig:geom_example}. The two invariants defining the space in this example are the \emph{chromatic number} (the minimum number of colors needed to color the vertices such that no two adjacent vertices share the same color) for the $X$-axis and the \emph{clique number} (the number of vertices in the largest complete subgraph) for the $Y$-axis. The points plotted represent all graphs of order 7. To illustrate that a single point can represent several non-isomorphic graphs, consider the coordinates in this figure: the point $(1, 1)$ corresponds uniquely to the empty graph on 7 vertices, and the point $(7, 7)$ corresponds uniquely to the complete graph on 7 vertices. In contrast, the point $(2, 2)$ represents 87 non-isomorphic graphs of order 7, a set that notably includes all 11 non-isomorphic trees on 7 vertices. Finally, this polytope makes an important inequality visually immediate: the chromatic number $\chi(G)$ of a graph $G$ is always greater than or equal to its clique number $\omega(G)$, that is the well-known relation $\chi(G) \ge \omega(G)$. 

Beyond observing inequalities for specific orders, the ultimate goal of this geometrical approach is to obtain a complete, parameterized polyhedral description valid for any graph order $n$. When such a general description can be formulated and mathematically proven, it yields the optimal set of linear inequalities governing the relationship between the two invariants. Indeed, by the very definition of a convex hull, its facets constitute the minimal set of bounding linear inequalities required to fully describe that relationship. This complete characterization has been successfully achieved for several pairs of invariants. For instance, it was established for the \emph{diameter} (the maximum eccentricity among all vertices) and the size of the graph in~\cite{melot2008facet}, as well as for the Fibonacci index (the total number of stable sets, where a stable set is a subset of pairwise non-adjacent vertices) and the \emph{independence number} (or stability number, that is the size of the largest stable set) in both general and connected graphs in~\cite{bruyere2009fibonacci}. Similarly, the complete polyhedral description for the independence number and the size of connected graphs was proven in~\cite{christophe2008linear}. Notably, this latter result formally resolved an open problem posed by Ore in 1962~\cite{ore1962theory}, which serves as a connected variant of the classical Tur\'an's theorem~\cite{turan1941extremalaufgabe}.

\section{Web interface and API}
\label{sec:interface_api}

This section details the current architecture and features of the PHOEG platform. Section~\ref{subsec:technologies} describes the underlying software stack. Section~\ref{subsec:functionalities} provides a comprehensive overview of the web interface functionalities, while Section~\ref{subsec:phoeg_api} presents the standalone API available for direct database access.

\subsection{Technologies}
\label{subsec:technologies}

While our research team has utilized the core PHOEG database for years, we recognized the need to make this resource accessible to the broader community, especially to researchers unfamiliar with formulating complex SQL queries. This objective drove the development of a dedicated web interface.

The platform follows a modern, three-tier architecture to ensure a clear separation of concerns: a frontend for user interaction, a backend Application Programming Interface (API) for business logic and query processing, and a relational database for reliable storage.

The frontend is built in TypeScript with React \cite{react}, leveraging its component-based architecture for efficient state management and dynamic UI updates. To ensure a clean, responsive, and accessible design, the interface heavily relies on Chakra UI~\cite{chakraui}, a robust component library.

The backend consists of a RESTful API developed in Rust~\cite{rust} using the \emph{rocket} framework~\cite{rocket}. This API acts as an intermediary, translating user interactions from the frontend into optimized database queries. Crucially, this API is also exposed for programmatic access independently of the web interface, as further detailed in Section~\ref{subsec:phoeg_api}.

The data layer is powered by a PostgreSQL relational database~\cite{postgresql}. It stores the invariant values for about 34 million pairwise non-isomorphic graphs including all graphs up to order $10$, along with other sparser families of graphs such as trees or chemical graphs (note that orders higher than 10 are not yet available within the interface). To guarantee that each graph is represented uniquely, the database relies on canonical forms computed via McKay's Nauty~\cite{mckay2014practical}. Currently, the database encompasses approximately fifty numeric invariants and a dozen Boolean ones. We continuously expand this collection as new research needs arise, and we actively encourage users to suggest the addition of new, relevant invariants.

\subsection{Functionalities of the website}
\label{subsec:functionalities}

This section details the core features of the web interface, publicly accessible at \url{https://phoeg.umons.ac.be}. The user workflow is logically divided into two main stages: defining the invariant space (Section~\ref{subsubsec:problem_definition}) and visually exploring the resulting polytopes and graphs (Section~\ref{subsubsec:problem_display}).

\subsubsection{Problem definition}
\label{subsubsec:problem_definition}

\begin{figure}[!htb]
\begin{center}
\includegraphics[width=0.6\textwidth]{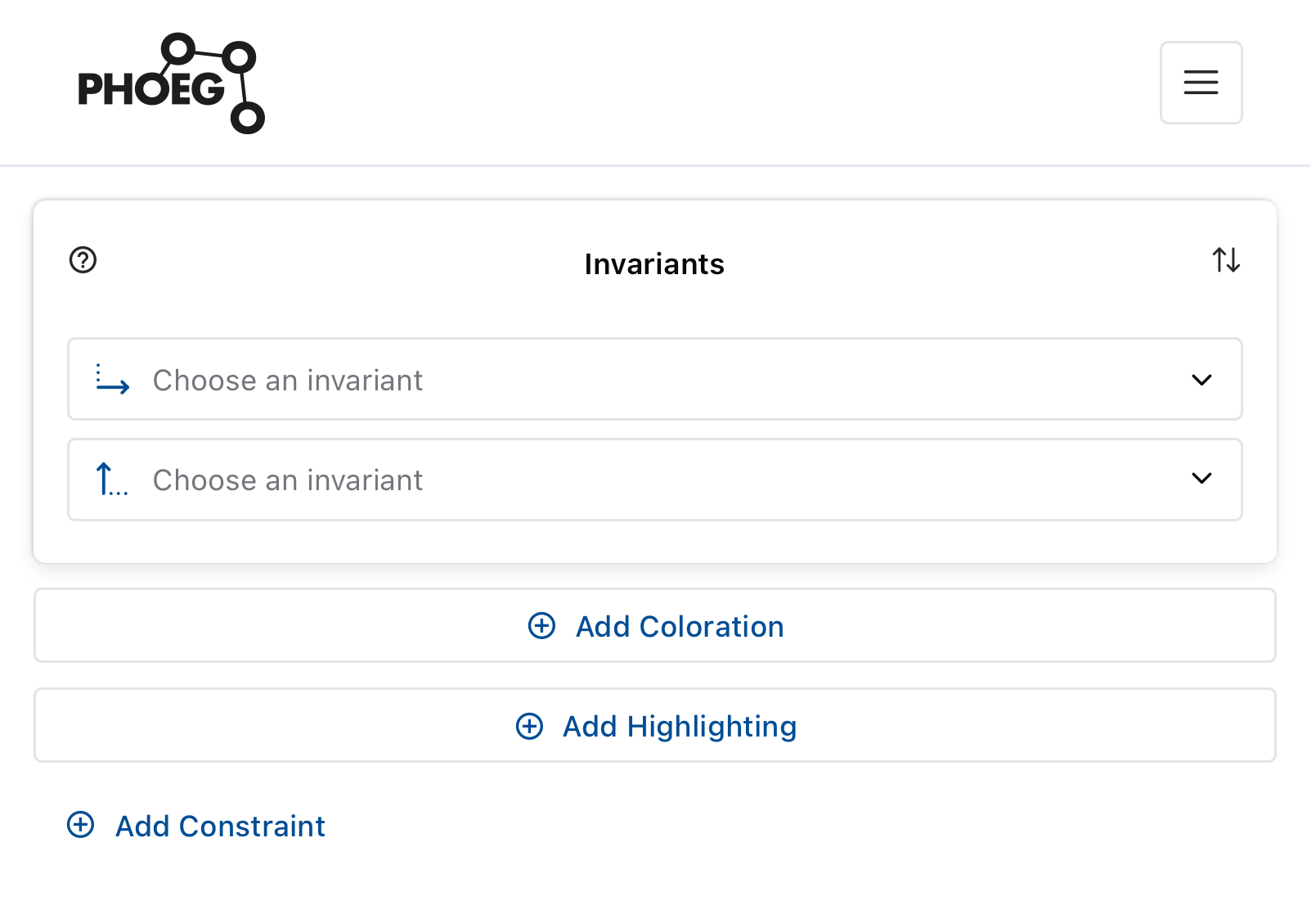}
\caption{Basic problem definition}\label{fig:problem_definition_simple}
\end{center}
\end{figure}

\begin{figure}[!htb]
    \centering

    \begin{subfigure}[c]{0.46\textwidth}
        \centering
        \includegraphics[width=\linewidth]{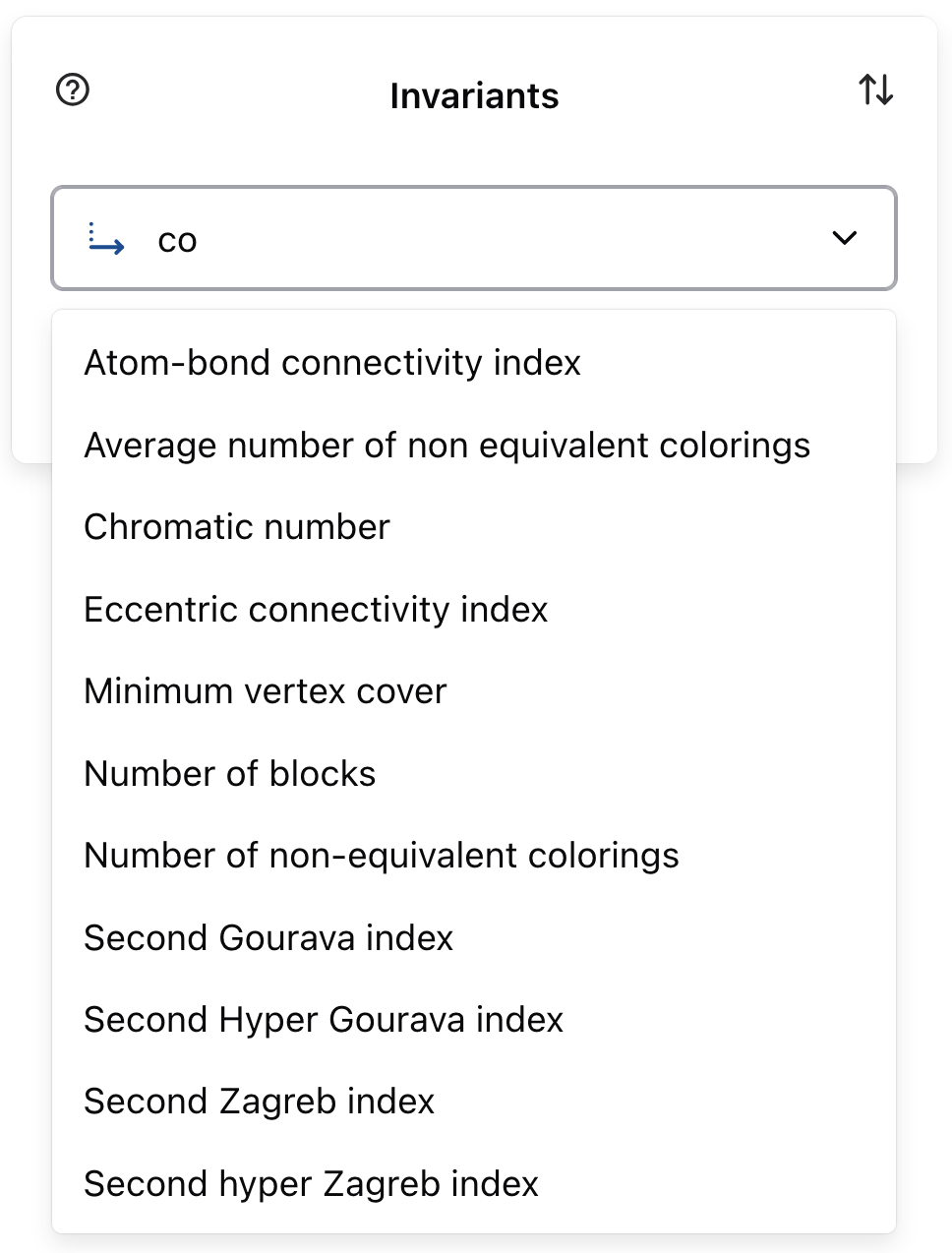}
        \caption{Invariant input field for the X-axis}
        \label{fig:invariant_input}
    \end{subfigure}
    \hfill
    \begin{minipage}[c]{0.46\textwidth}
        
        \begin{subfigure}{\linewidth}
            \centering
            \includegraphics[width=\linewidth]{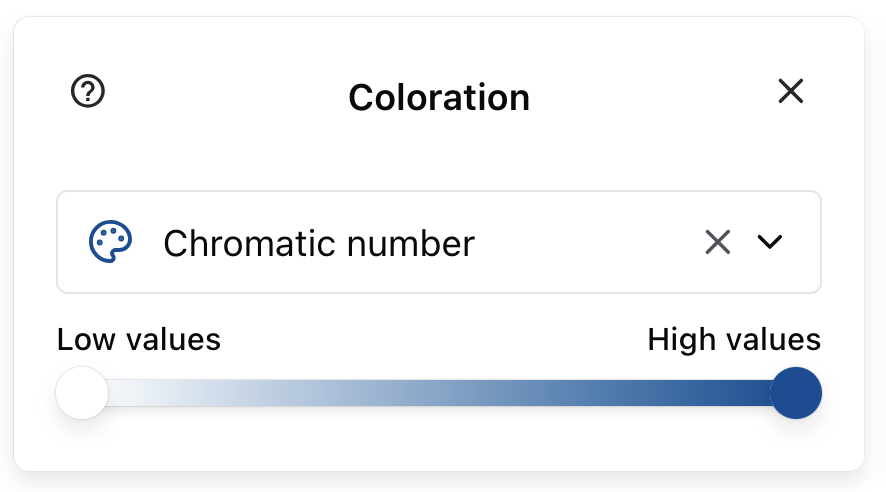}
            \caption{Defining the coloration}
            \label{fig:coloration_definition}
        \end{subfigure}
        
        \vspace{0.6cm}
        
        \begin{subfigure}{\linewidth}
            \centering
            \includegraphics[width=\linewidth]{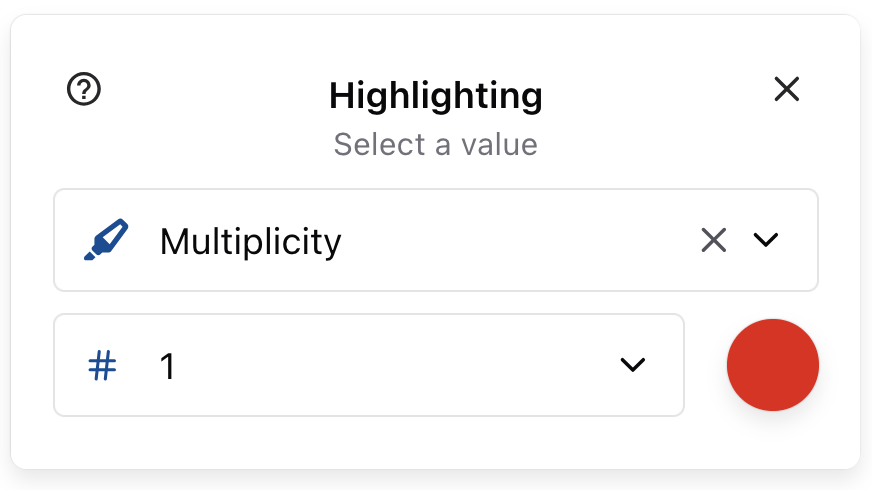}
            \caption{Defining the highlighting}
            \label{fig:highlighting_definition}
        \end{subfigure}
        
    \end{minipage}
   
    \caption{Overview of the parameter input boxes in the Problem Definition panel.}
    \label{fig:problem_definition_inputs}
\end{figure}

Upon accessing the homepage, the user's first task is to configure the mathematical problem, as shown in Figure~\ref{fig:problem_definition_simple}. The only mandatory step is defining the two-dimensional invariant space by selecting the invariants for the $X$ and $Y$ axes. The input fields feature auto-completion to assist the user in navigating the extensive list of available invariants (see Figure~\ref{fig:invariant_input}). Once these axes are defined, the resulting polytope is dynamically generated and displayed (see Section~\ref{subsubsec:problem_display}).

To refine the exploration, the user can incrementally add three types of optional parameters:
\begin{itemize}
    \item \textbf{Coloration:} Applies a continuous color gradient to the polytope points based on the value of a third selected invariant. The user defines a bi-color gradient representing the minimum and maximum values of this invariant (Figure~\ref{fig:coloration_definition}). 

    \item \textbf{Highlighting:} Isolates and distinctly colors points (and their corresponding drawn graphs) that strictly match a specific target value for a chosen invariant (see for example Figure~\ref{fig:highlighting_definition} where the points that represent graphs with an average number of non-equivalent colorings of $2$ (and drawn graphs meeting the same criterion) will be colored red).

    \item \textbf{Constraints:} Acts as a dataset filter, restricting the visualized graphs to those satisfying specific Boolean conditions (true/false) or numerical bounds (greater than, less than, or equal to a target value). Figures~\ref{fig:numerical_constraint} and~\ref{fig:boolean_constraint} illustrate these constraint definition boxes.

\begin{figure}[!htb]
    \centering
    \begin{subfigure}[b]{0.46\textwidth}
        \centering
        \includegraphics[width=\linewidth]{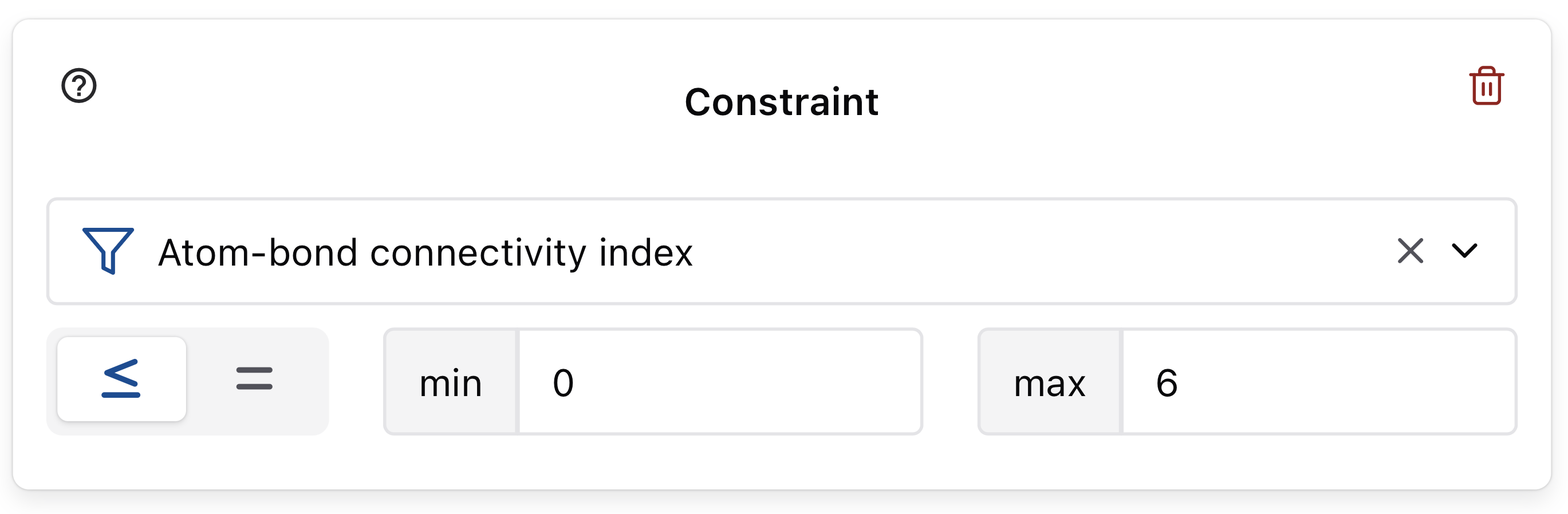}
        \caption{Numerical constraint}
        \label{fig:numerical_constraint}
    \end{subfigure}
    \hfill
    \begin{subfigure}[b]{0.46\textwidth}
        \centering
        \includegraphics[width=\linewidth]{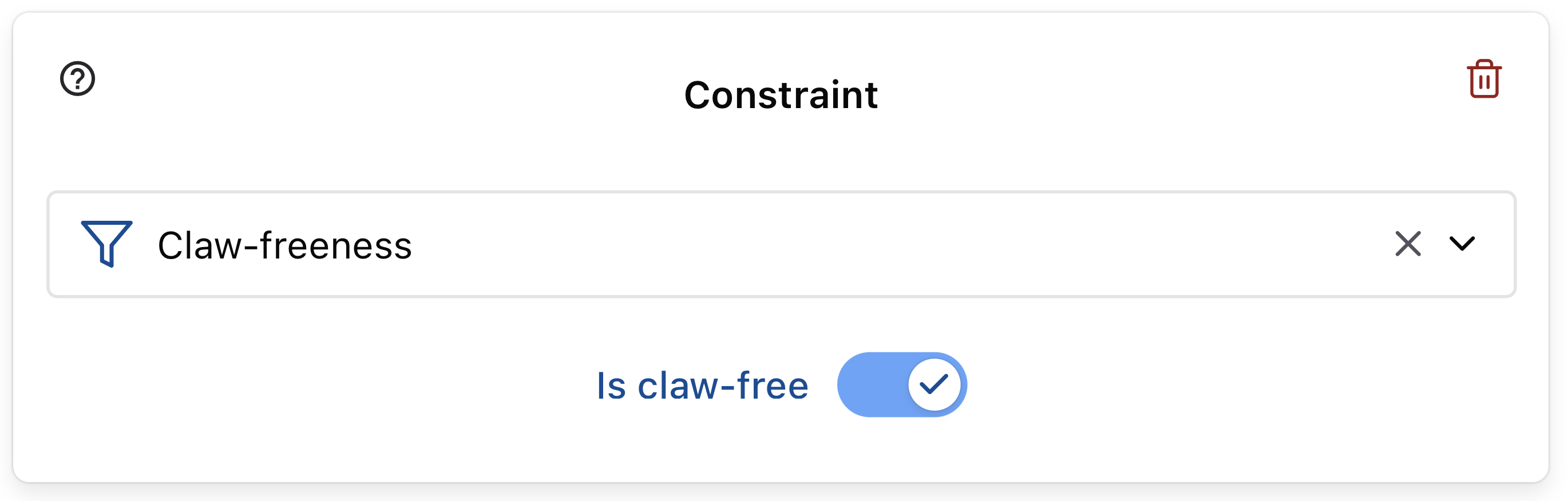}
        \caption{Boolean constraint}
        \label{fig:boolean_constraint}
    \end{subfigure}

    \caption{Examples of constraint definition boxes.}
    \label{fig:constraints_boxes}
\end{figure}
\end{itemize}

Examples of practical applications of coloration and highlighting are provided in the next section. Finally, it is worth noting a special option available within the coloration and highlighting menu: \emph{Multiplicity}. While not a mathematical graph invariant per se, multiplicity acts as a powerful visual tool by mapping the color scale (or highlight) to the density of the points---specifically, the exact number of pairwise non-isomorphic graphs residing at a single $(x, y)$ coordinate. This allows users to instantly distinguish visually between rare (e.g., points with a multiplicity of 1) and dense hubs where numerous graphs overlap.

\subsubsection{Problem display}
\label{subsubsec:problem_display}

\begin{figure}[!ht]
\begin{center}
\includegraphics[width=0.9\textwidth]{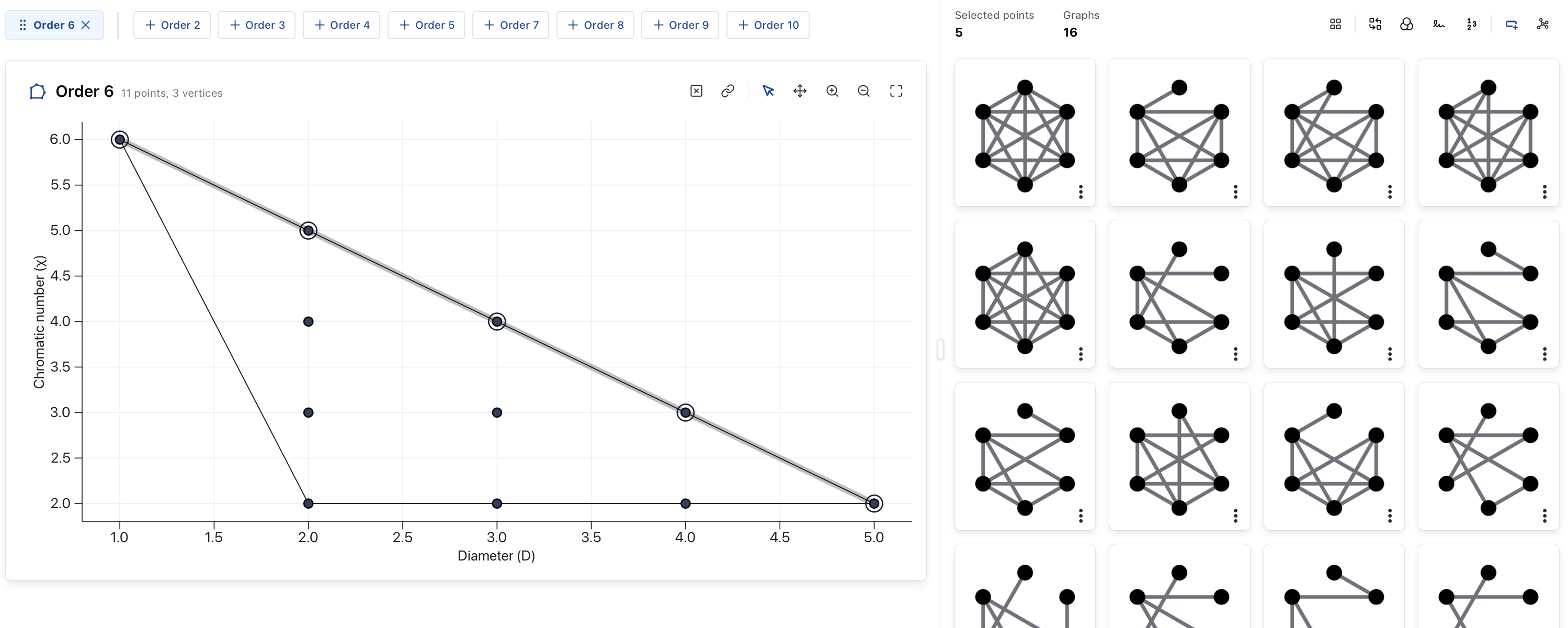}
\caption{Problem display interface}\label{fig:problem display}
\end{center}
\end{figure}

\begin{figure}[!ht]
\begin{center}
\includegraphics[width=0.9\textwidth]{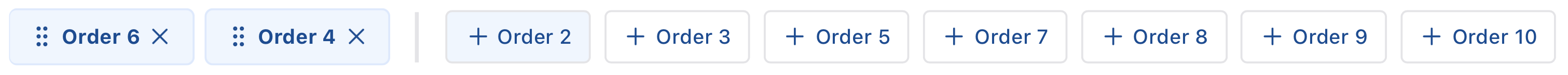}
\caption{Order selection}\label{fig:order_selection}
\end{center}
\end{figure}

Once the problem is defined, the exploration interface is populated (Figure~\ref{fig:problem display}). The user can filter and reorder the visualization of the polytopes by associated graph orders using toggle buttons at the top left (Figure~\ref{fig:order_selection}). Polytopes from order 2 to 10 are available, with order 6 displayed by default.

The interface is divided into two highly interactive panels. To ensure optimal usability across any screen size, users can manually adjust the relative width of the left and right panels, with all resulting polytopes and graph drawings scaling automatically to fit the available space:
\paragraph{Polytope Exploration (Left Panel):} This section displays and allows to scroll the polytopes for the selected orders, constructed from the appoach described in Section~\ref{sec:geometrical_approach}. The axes reflect the chosen invariants, and metadata (such as the number of points and vertices) is provided. Users can navigate the space by zooming, panning, or synchronizing axes across different orders for direct scale comparison. Interactivity is a core feature: clicking a specific point selects it, while clicking a facet automatically selects all its incident points. 

\begin{figure}[!ht]
\begin{center}
\includegraphics[width=0.7\textwidth]{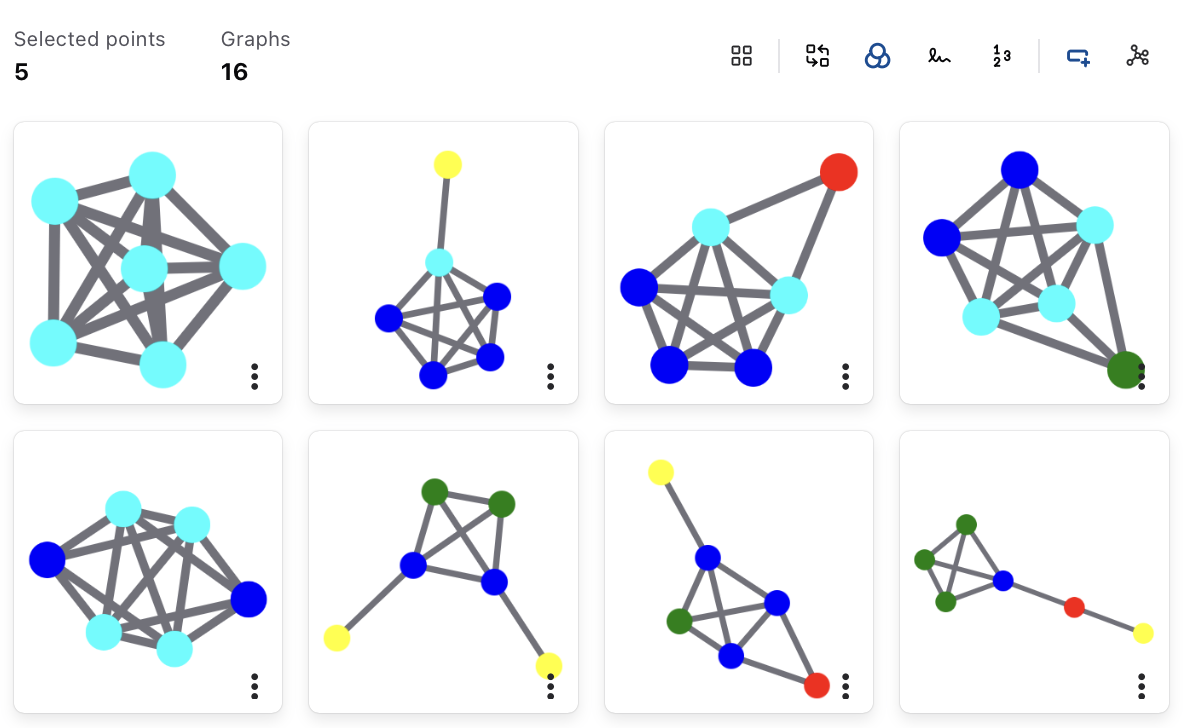}
\caption{Layout using the $cose$ algorithm and with a degree-based coloration of the vertices}\label{fig:drawn_graphs_cose}
\end{center}
\end{figure}

\begin{figure}[!ht]
\begin{center}
\includegraphics[width=0.6\textwidth]{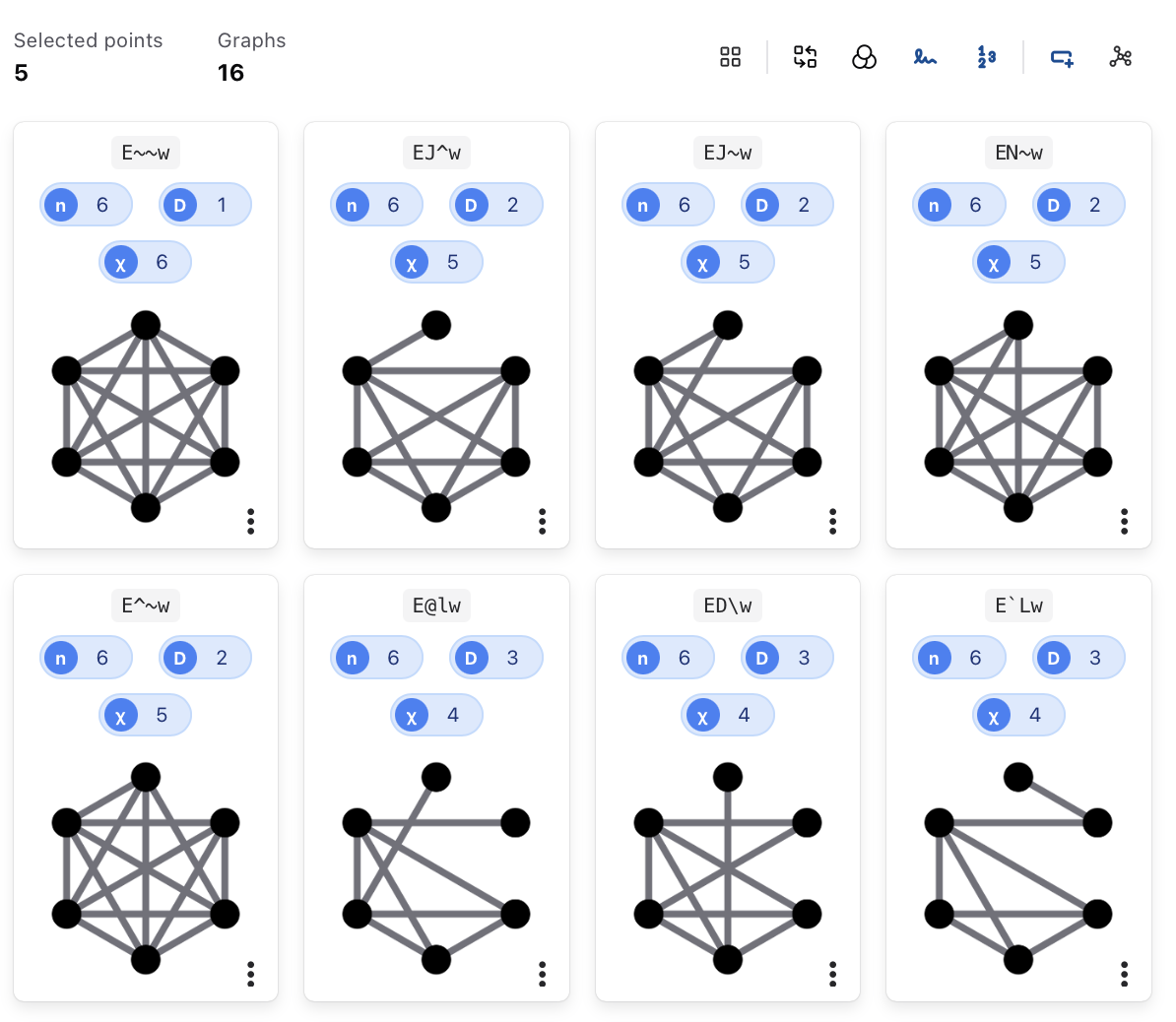}
\caption{Graph rendering with signature and invariant annotations}\label{fig:drawn_graphs_more}
\end{center}
\end{figure}

\begin{figure}[!ht]
\begin{center}
\includegraphics[width=0.7\textwidth]{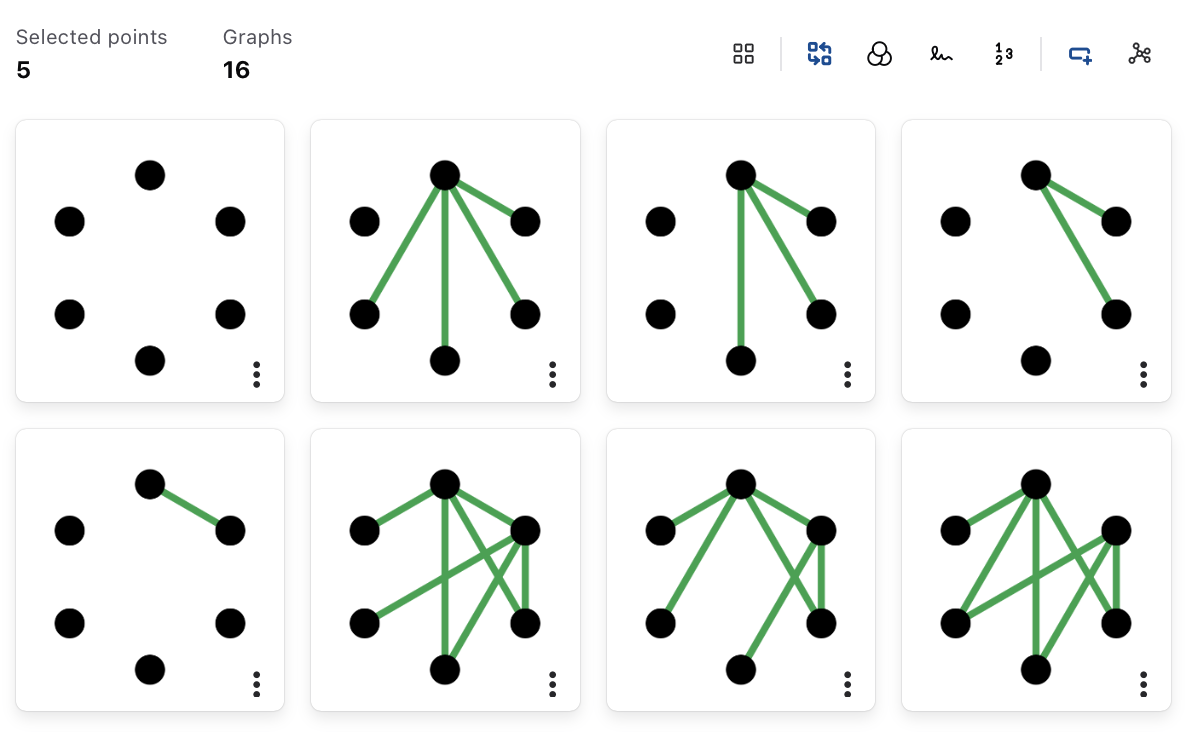}
\caption{Complement graphs rendering}\label{fig:drawn_graphs_complement}
\end{center}
\end{figure}

\paragraph{Graph Visualization (Right Panel):} This section dynamically renders the graphs corresponding to the points selected in the left panel. To handle large datasets, graphs are lazy-loaded upon scrolling. The graph drawing module offers extensive customization:
\begin{itemize}
    \item \textbf{Layout and Rendering:} Users can drag vertices manually or apply one of 8 automated layout algorithms (e.g., \emph{circle} as seen in Figure~\ref{fig:problem display} or \emph{cose}, as seen in Figure~\ref{fig:drawn_graphs_cose}). If necessary, the user can drag vertices around to move them.
    \item \textbf{Graph Data:} By default, only graph drawings are shown, but the user can toggle the display to reveal graph signatures, invariant values (Figure~\ref{fig:drawn_graphs_more}), or render complement graphs (where edges and non-edges are inverted and displayed in green, see Figure~\ref{fig:drawn_graphs_complement}).
    \item \textbf{Visual Cues:} Vertex colors can be mapped to their respective degrees for immediate structural insight. For example, the graphs  drawn in Figure~\ref{fig:drawn_graphs_cose} have their vertices colored depending on their degree (degree $1$ in yellow, degree $2$ in red, etc.).
\end{itemize}
These visual settings can be applied globally to all displayed graphs or individually via a context menu below each graph.

\begin{figure}[!ht]
\begin{center}
\includegraphics[width=0.9\textwidth]{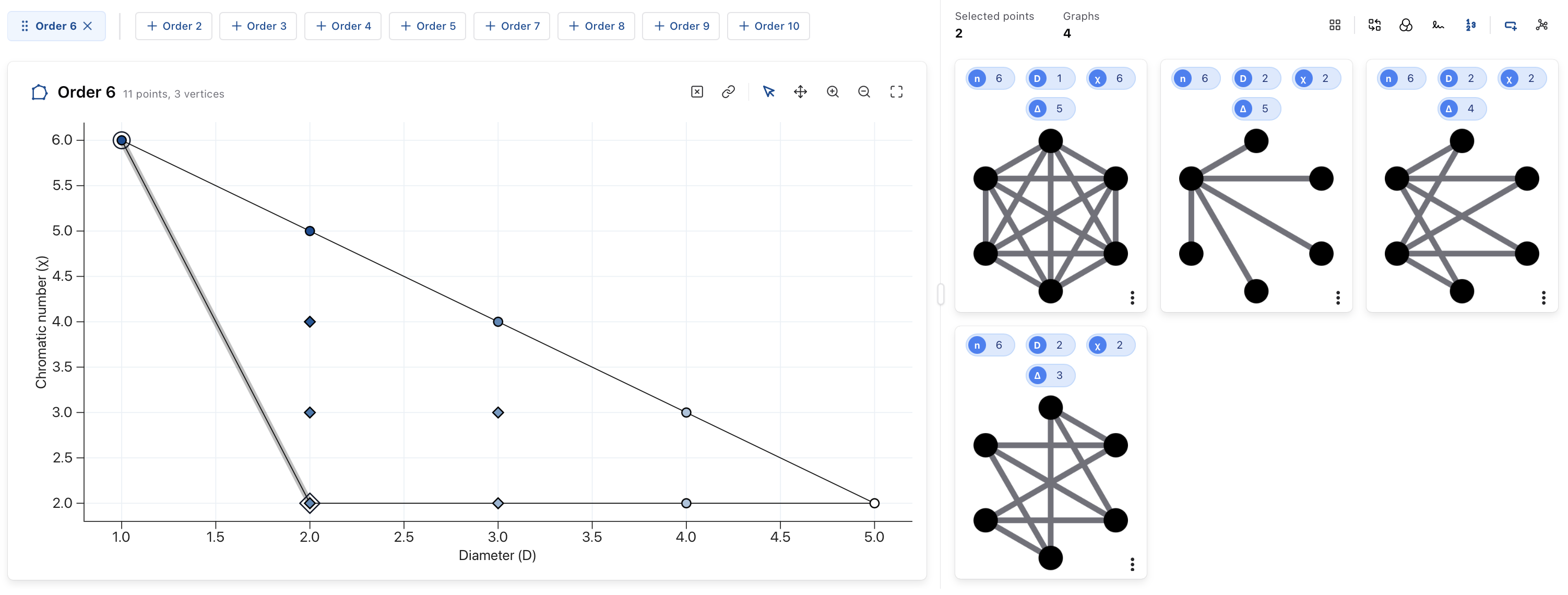}
\caption{Coloration based on the maximum degree, with a gradient from white to blue}\label{fig:display_coloration}
\end{center}
\end{figure}

\begin{figure}[!ht]
\begin{center}
\includegraphics[width=0.9\textwidth]{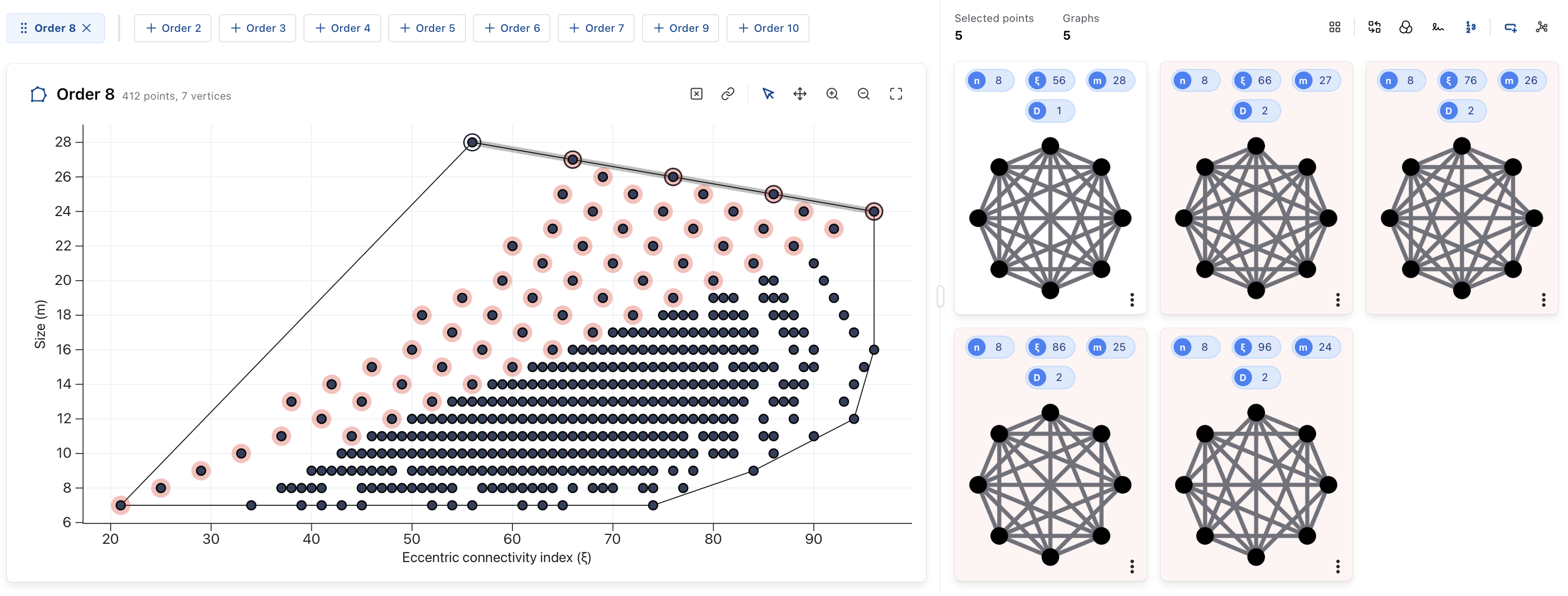}
\caption{Highlighting graphs with a diameter of $2$}\label{fig:display_highlighting}
\end{center}
\end{figure}

Furthermore, if a \emph{Coloration} or \emph{Highlighting} was defined in the previous step, it is visually applied in both panels. For instance, Figure~\ref{fig:display_coloration} illustrates a coloration based on the \emph{maximum degree} (the highest number of edges incident to any single vertex in the graph), applying a gradient from white to blue across the points of the polytope. As established in Section~\ref{sec:geometrical_approach}, a single geometric coordinate can represent multiple non-isomorphic graphs. If these underlying graphs yield different values for the invariant used for coloration, the tool cannot assign a single color; instead, it renders the point as a diamond to alert the user of this variance, a feature clearly visible in the same figure. To properly understand Figure~\ref{fig:display_highlighting}, we first introduce two notions. The \emph{eccentricity} $\varepsilon(v)$ of a vertex $v$ is defined as the greatest distance between $v$ and any other vertex in the graph. The \emph{Eccentric Connectivity Index} (ECI) of a graph $G$, denoted by $\xi^c(G)$, is the sum over all vertices of the product of their degree and their eccentricity, that is, $\xi^c(G) = \sum_{v \in V(G)} d(v)\,\varepsilon(v)$. Figure~\ref{fig:display_highlighting} demonstrates a highlighting rule applied to the polytope of graphs of order 8, plotting the eccentric connectivity index against the size. Specifically, points corresponding to graphs with a diameter of 2 are colored red. This feature allows to immediately identify that graphs of diameter 2 are clustered in the upper region of the polytope. The highlighting is also applied to the right panel, where the corresponding graph drawings are displayed with a red background. In Figure~\ref{fig:display_highlighting}, the complete graph on 8 vertices is the only non-highlighted graph since its diameter is 1.

Finally, a notable feature of the web interface is that most of the session state is encoded in the URL. This encompasses not only the full problem definition (selected invariants, applied constraints, coloration, and highlighting parameters) but also most of the current state of the display: the selected graph orders, the points selected in the polytopes, whether axis synchronization is enabled, and the active graph visualization options. This has several practical implications: the interface is inherently refresh-resistant, as reloading the page faithfully restores the exact state of the session. It also makes problems trivially shareable. To streamline this process, the interface features a dedicated \textbf{Share Configuration} button located in the problem definition panel. Clicking this button automatically copies the complete, state-encoded URL to the user's clipboard. Anyone opening this link is brought to the exact same analytical view, making it an ideal tool for collaborative research, saving progress, or linking directly to a specific dataset from an academic paper. For instance, if you follow \href{https://phoeg.umons.ac.be/phoeg/?q=N4IgHiBcIKYMYEsQBoQE8ogHYFcC2A%2BjACYDmMAziiHFANoC6qAFlKEtMQgIZ7W3QAxDADsARgAcABgBG1AG5QATAF9UcADZsQWoQDMDh6qyFSpAFmIBOAKwgVKoA&s=N4Ig9iBcDaAcC6AaEBnK1oE4AsiCMsArIgEyYDMi5sADItnnvEiAGbrShqR7ICmUEgF9myAOZRQAEygAXAE4BXPsgA2UEAGMAlvM2qByWehAA7RQFsA%2BgDc%2B82ds180-HSGTnrfKWJceQKW0AQwsQeCEhIA}{this link}, you will get the exact configuration used in Figure~\ref{fig:display_highlighting}.

\subsection{API}
\label{subsec:phoeg_api}

Beyond the web interface, PHOEG exposes its core functionalities through a standalone RESTful API. This programmatic access is designed for researchers who wish to integrate PHOEG's extensive database into their own computational workflows, scripts, or automated pipelines without manual interaction with the UI.

The API serves as a bridge to the PostgreSQL database, allowing users to perform complex queries using standard HTTP requests. Data is typically exchanged in JSON format, ensuring compatibility with most modern programming languages such as Python or R. Key capabilities of the API include:
\begin{itemize}
    \item \textbf{Data Retrieval:} Fetching specific invariant values for a given graph identified by its canonical signature.
    \item \textbf{Advanced Filtering:}
    Retrieving the individual graphs lying on specific points of the polytope. Given a pair of invariants defining the two axes, a graph order, and a list of $(x,y)$ coordinates, it returns all graphs of that order whose invariant values match those coordinates. The result set can be further narrowed using range constraints on any invariant, and optionally enriched with values of additional invariants. A coloration invariant can also be specified, whose value is then attached to each returned graph to support the coloration feature of the interface.
    \item \textbf{Exporting:} Obtaining graph data in the \emph{graph6} standard format~\cite{graph6} for further local analysis or for use with other graph theory software.
\end{itemize}

The current API is documented at \url{https://phoeg.umons.ac.be/phoeg/api/docs}.

\section{Using PHOEG in research}
\label{sec:phoeg_research}

PHOEG can be utilized in research across various contexts. It serves as an excellent starting point for the preliminary inspection of a problem, helping researchers to build intuition. Crucially, it assists in formulating new conjectures, as detailed in Section~\ref{subsec:find_conjectures}, in discovering counterexamples to existing ones, as demonstrated in Section~\ref{subsec:find_counterexamples}, and even provides assistance for formal proofs, as discussed in Section~\ref{subsec:help_for_proofs}.

\subsection{Finding conjectures}
\label{subsec:find_conjectures}

The visual and geometric nature of PHOEG acts as a powerful catalyst for forming new mathematical hypotheses. By visualizing the invariant space and explicitly identifying the graphs located on the bounding facets, researchers can gather critical insights even when restricted to graphs of small orders. This empirical observation often makes it surprisingly straightforward to deduce the structural properties of extremal graphs for any arbitrary order $n$, and subsequently, to generalize the corresponding inequalities. In practice, when computing and displaying polytopes for successive graph orders, it is remarkably common to observe that their overall geometric shapes remain structurally consistent. These recurring polyhedral patterns immediately capture the observer's attention, significantly easing the extrapolation of vertex coordinates and the formulation of general, parameterized conjectures.

\begin{figure}[!ht]
\begin{center}
\includegraphics[width=0.32\textwidth]{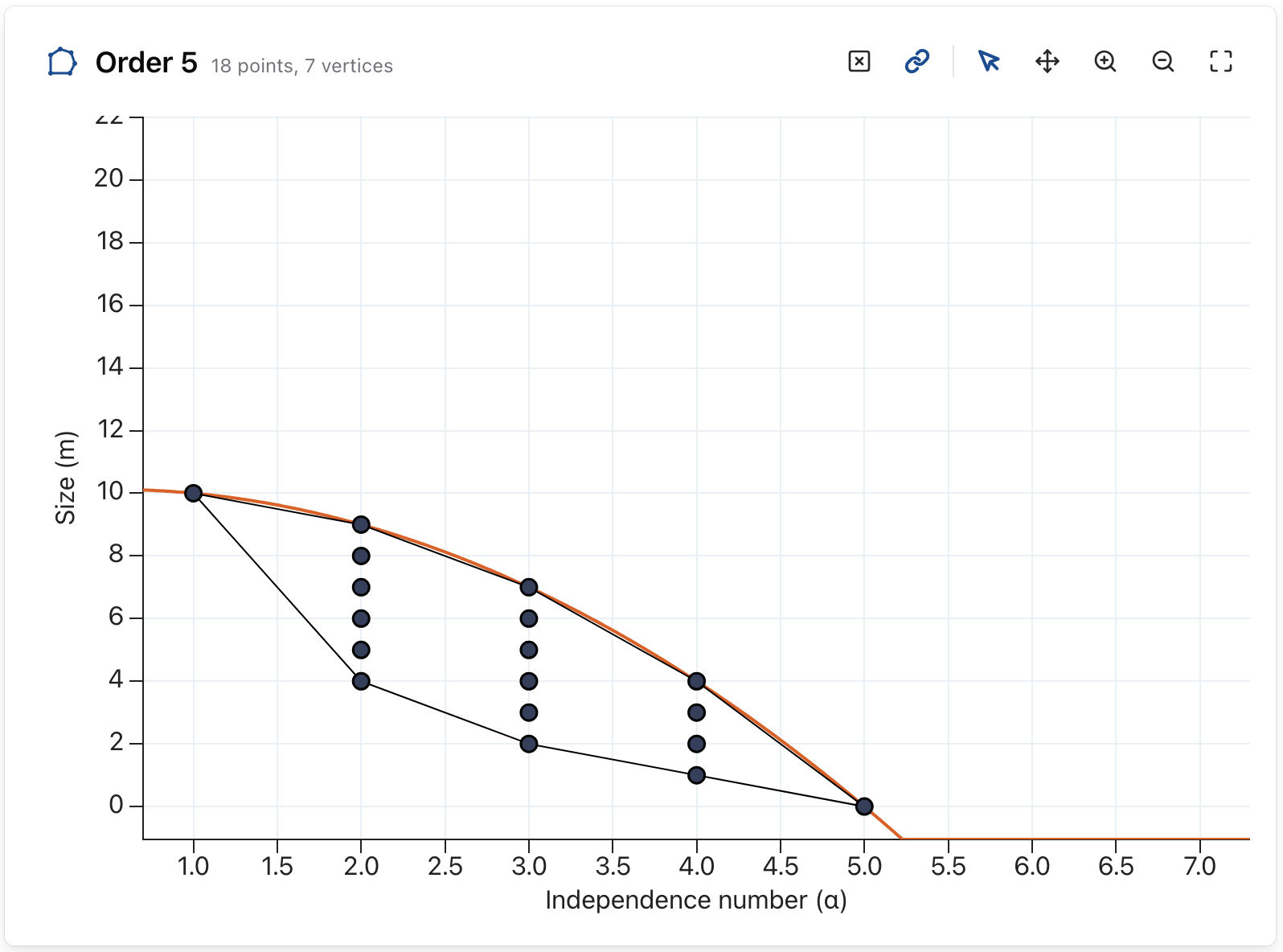} \ \includegraphics[width=0.32\textwidth]{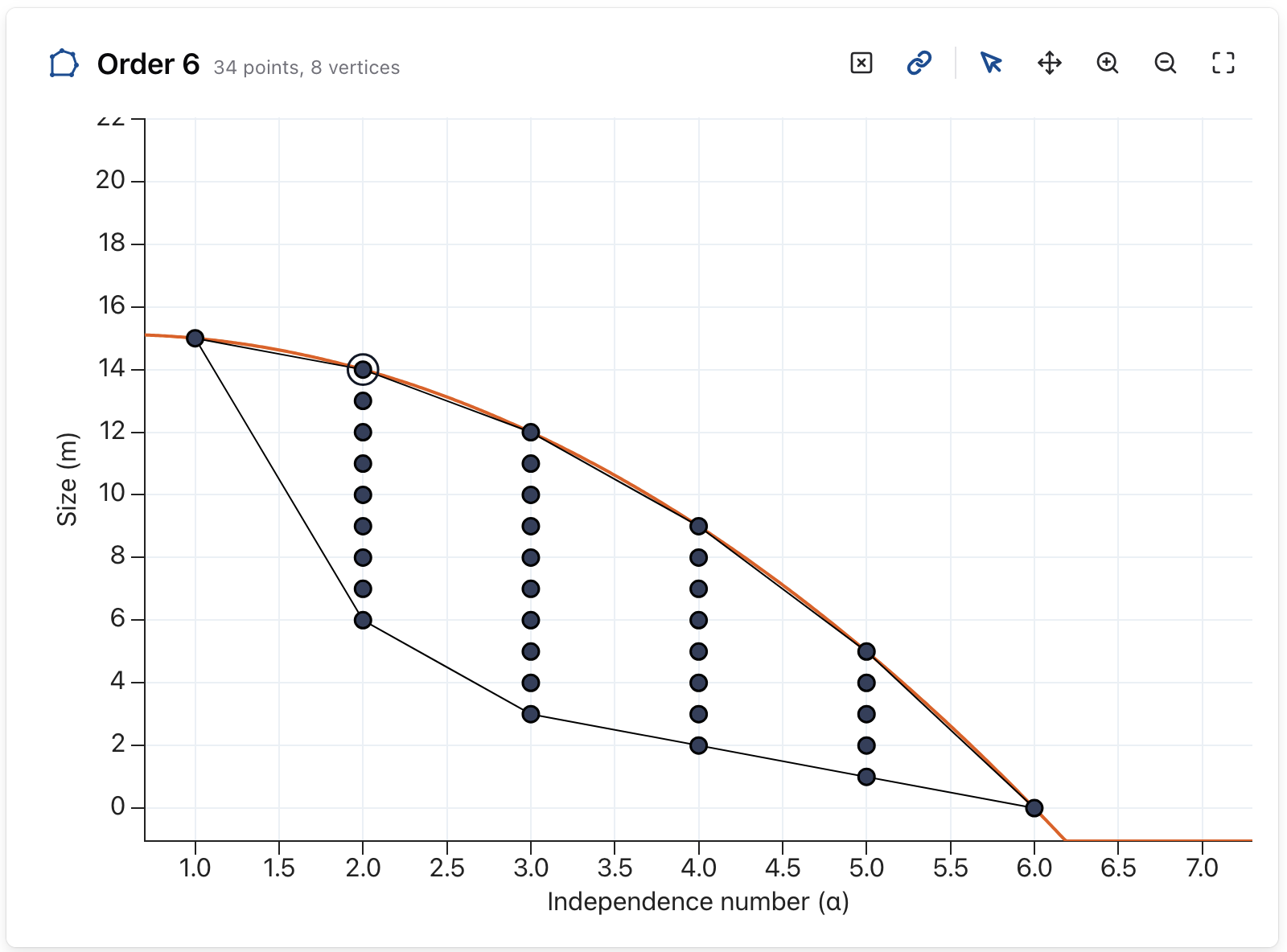} \ \includegraphics[width=0.32\textwidth]{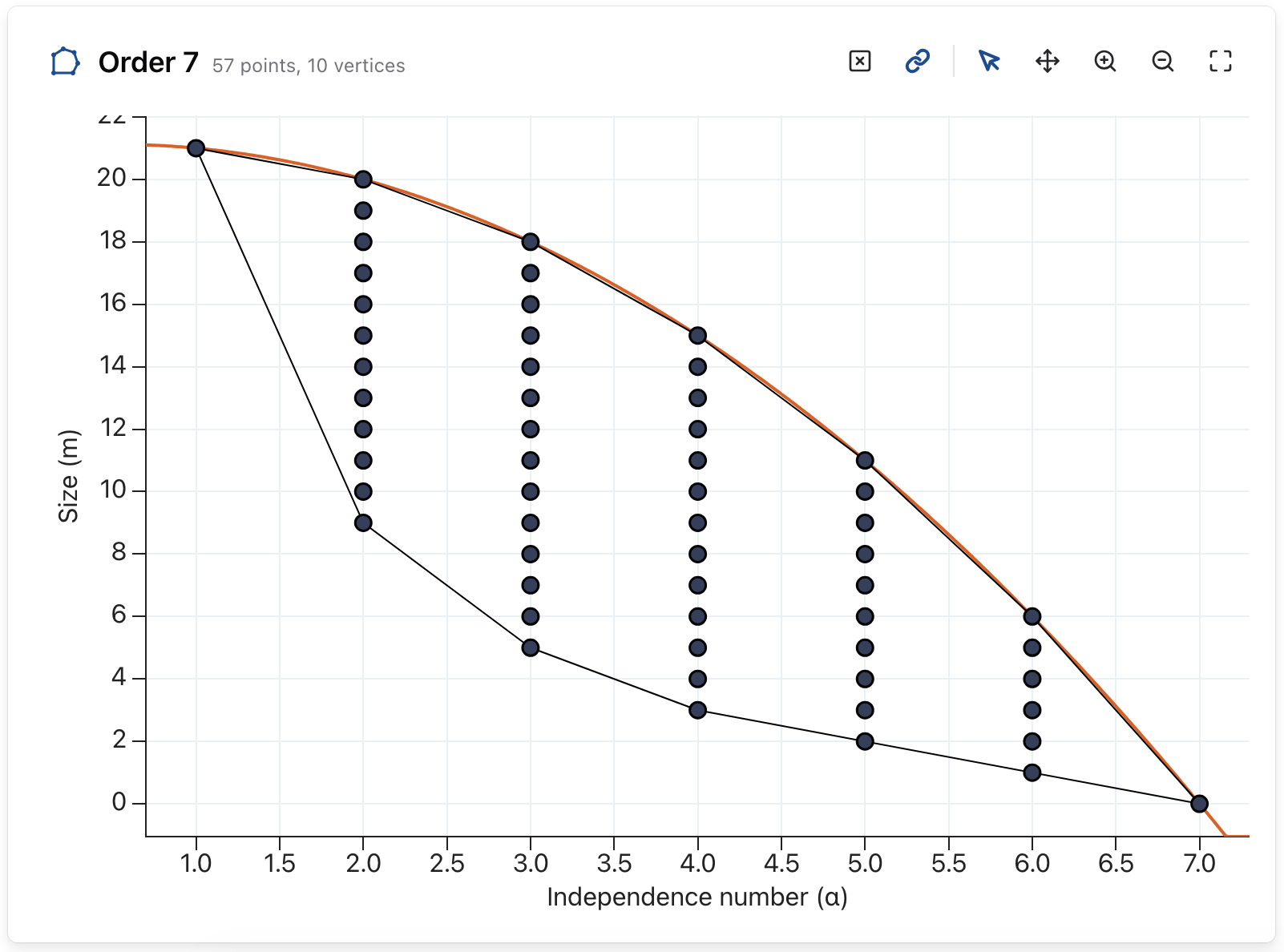} 
\caption{Independence number against the number of edges for graphs of order 6, 7, and 8.}\label{fig:turan_polytopes}
\end{center}
\end{figure}

An illustration of this phenomenon is provided in Figure~\ref{fig:turan_polytopes}, which displays the polytopes for graphs of orders 6, 7, and 8, mapping the independence number $\alpha$ against the size $m$. The consistent geometric shape across these successive orders is immediately apparent. By interacting with the lower boundary of these polytopes—representing the problem of minimizing the number of edges for a given independence number—a user can easily inspect the corresponding extremal graphs. Such an examination quickly reveals that these vertices consistently correspond to Tur\'an graphs (formed by disjoint unions of cliques of balanced sizes), perfectly aligning with the well-known result established by Tur\'an~\cite{turan1941extremalaufgabe}. 

To further support hypothesis formulation and testing, PHOEG includes a \emph{Plot Function} feature. This tool allows users to superimpose custom algebraic curves onto the two-dimensional space by defining a formula that relates the $X$- and $Y$-axis invariants, optionally incorporating the graph order $n$. As an illustration, while the lower boundary of the polytopes in Figure~\ref{fig:turan_polytopes} is populated by Tur\'an graphs, the upper boundary represents the maximization of the size $m$ for a fixed independence number $\alpha$. Using this plotting feature, we overlaid the function $m = \frac{1}{2} \left(n^2 - n - \alpha^2 + \alpha \right)$ (in red in Figure~\ref{fig:turan_polytopes}). This curve perfectly aligns with the upper extremal vertices, which correspond to \emph{complete split graphs}---graphs whose vertex set can be partitioned into a clique and an independent set, with every vertex in the clique adjacent to every vertex in the independent set.

This core methodology has a long track record of driving research. Its predecessor, GraPHedron, was used internally by our research team starting in 2005, before being progressively replaced by PHOEG in the mid-2010s. When the geometrically derived conjectures are successfully proven, they contribute directly to the broader field of extremal graph theory. Specifically, the hypotheses generated and refined through GraPHedron and PHOEG have led to several peer-reviewed publications. Notable examples include results on the independence number~\cite{christophe2008linear, bruyere2009fibonacci, bruyere2012trees}, the number of non-equivalent colorings~\cite{hertz2016counting, absil2018sharp}, the Eccentric Connectivity Index~\cite{devillez2019minimum, hauweele2019maximum}, an analysis of the approximation factor of a maximal matching heuristic~\cite{cardinal2005tight}, extremal properties in chemical graph theory~\cite{hertz2025extremal, bonte2026extremal}, as well as invariants concerning the average number of colorings and the average number of matchings~\cite{hertz2024average, hertz2023lower, hertz2023upper}. A comprehensive and continuously updated list of all related publications is maintained on the PHOEG website.

\subsection{Finding counterexamples}
\label{subsec:find_counterexamples}

Beyond original discovery, PHOEG is a highly effective platform for verifying existing mathematical claims. If a known conjecture is formulated in a way that is compatible with our two-dimensional geometric approach—typically involving bounds or relationships between invariants—it can be systematically tested within the interface. By generating the relevant polytope, researchers can quickly gather empirical support for the conjecture, refine its proposed bounds, or outright refute it by uncovering unexpected points that violate the hypothesized inequalities. 

To give a concrete example, in 2014, Zhang et al.~\cite{zhang2014maximal} proposed a conjecture regarding the Eccentric Connectivity Index (ECI). In this context, for given values $(n, m)$, they defined $$d = \left\lfloor \dfrac{2n+1-\sqrt{17+8(m-n)}}{2} \right\rfloor$$ and $E_{n,m}$ as the graph obtained from a clique $K_{n-d-1}$ and a path 
$
P_{d+1} = v_0 v_1 \ldots v_d
$
by joining each vertex of the clique to both $v_d$ and $v_{d-1}$, and by joining 
$
m - n + 1 - {n - d \choose 2}
$
vertices of the clique to $v_{d-2}$. The conjecture is as follows:

\begin{conj}[Zhang, Liu and Zhou, \cite{zhang2014maximal}]
Let $d_{n,m} \geq 3$. Then $E_{n,m}$ is the unique graph with maximal eccentric connectivity index among all connected graphs with $n$ vertices and $m$ edges.
\end{conj}

\begin{figure}[!ht]
\begin{center}
\includegraphics[width=0.9\textwidth]{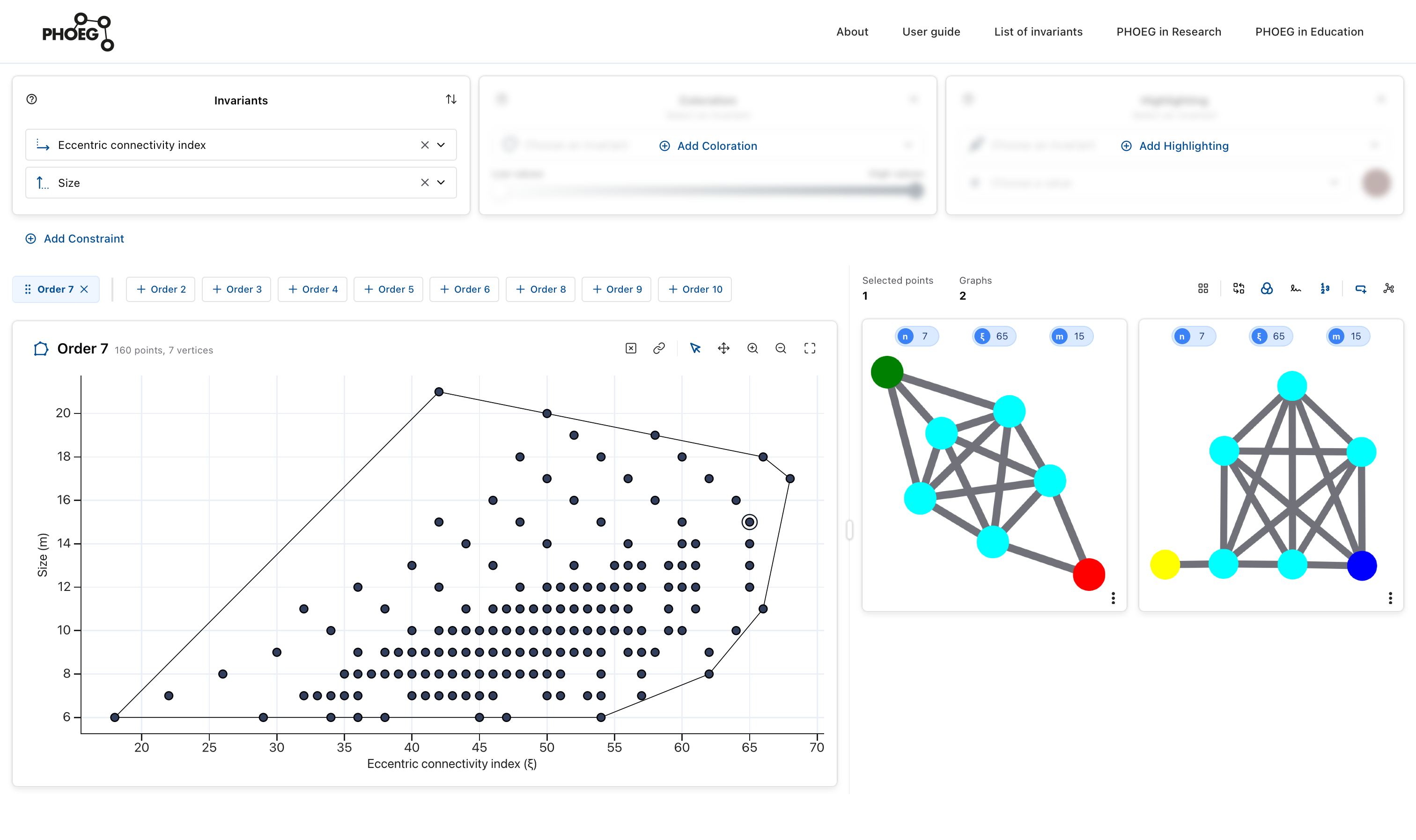}
\caption{Counterexample for Zhang's conjecture for $n=7, m=15, d=3$}\label{fig:counterexample_eci}
\end{center}
\end{figure}

The authors prove that the conjecture is true for certain values of $n$ and $m$. However, with the help of PHOEG, we can find a counterexample with little effort. Figure~\ref{fig:counterexample_eci} shows this counterexample. The selected point corresponds to $n=7$ and $m=15$, which gives $d_{7,15}=3$. This point corresponds to two graphs that have an ECI of $65$. The graph on the right is $E_{7,15}$, but the other graph is not isomorphic to $E_{7,15}$, which shows that $E_{7,15}$ is not the unique graph with maximal eccentric connectivity index among all connected graphs with $7$ vertices and $15$ edges.

\subsection{Help for proofs}
\label{subsec:help_for_proofs}

Beyond formulating conjectures, visually observing extremal points and their corresponding graph structures can provide critical intuition for formal demonstrations, particularly for proofs by transformation. In this mathematical approach, a graph is modified step-by-step to reach an extremal configuration while strictly increasing or decreasing a given invariant.

To support this methodology, a dedicated computational module named TransProof was developed as part of Devillez's PhD thesis~\cite{devillez2022proofs}. Although TransProof is not currently available within the public PHOEG web interface, its mechanics are deeply tied to the project's ecosystem. The primary function of TransProof is to compute a \emph{metagraph} for a specified set of graphs and structural transformations. In this metagraph, vertices represent the graphs themselves, and directed arcs represent the applied transformations.

Once computed, this metagraph is stored in a database, allowing researchers to run complex queries to test and refine their proof strategies. For instance, a user can query the database to check for counterexamples to a proposed proof idea by ensuring the absence of non-improving arcs. Furthermore, the flexibility of the database query language allows for more advanced topological checks, such as verifying whether every non-extremal graph possesses at least one strictly improving arc, or identifying if certain non-extremal graphs cannot be improved by the considered set of transformations.

\section{Using PHOEG in education}
\label{sec:phoeg_education}

Beyond its applications in research, PHOEG is actively utilized for educational purposes at both the University of Mons (UMONS) and Polytechnique Montr\'eal, specifically within the context of advanced graph theory courses. 

To facilitate the onboarding process for students and new researchers alike, the platform includes a comprehensive suite of built-in interactive tutorials. Accessible directly from the interface's help section, these step-by-step guides allow users to familiarize themselves with the tool at their own pace. Through hands-on, guided sandboxes, users learn the core geometric concepts, how to apply dataset constraints and colorations, and how to effectively navigate the polytope and graph panels before tackling more complex, unguided assignments.

For instance, at UMONS, in the \emph{Graphs \& Artificial Intelligence} course, once students have completed these introductory tutorials, the tool is further explored through a practical session where students are invited to use PHOEG's web interface to explore and solve classical problems in extremal graph theory. 

The first problem assigned to the students is based on a fictitious police investigation involving Al Capone, presented as a direct continuation of a problem introduced earlier in the course. In that earlier problem, an inspector had met seven mafiosi: John Gotti, Joe Adonis, Albert Anastasia, Liborio Bellomo, Tommaso Buscetta, Giuseppe Calicchio, and Al Capone. Each of them declared the number of people they had business relations with, except Al Capone, who remained deliberately vague and only admitted to having 'more than $3$' such relations. Using the handshaking lemma, students had previously established that the sum of all degrees must be even, which forces Al Capone to have exactly $5$ business relations, yielding the full degree sequence $2, 3, 3, 3, 4, 4, 5$.

The investigation has since progressed: the police have almost proven that Al Capone is guilty of murder, but he claims he was with Albert Anastasia at the time of the crime, providing himself an alibi. A new piece of evidence has also emerged: the largest group of people with no business relations between them is of size $2$. In graph-theoretical terms, this clue is meant to inform them that the independence number of the graph representing this situation is of size $2$. The objective of the riddle is to determine whether Al Capone (the vertex of degree $5$) could have a business relationship with Albert Anastasia (the vertex of degree $2$).

By applying these constraints in PHOEG, students discover that only a single graph matches the entire description. Because this unique graph contains exactly one vertex of degree $5$ and one vertex of degree $2$, they can directly inspect its structure to see if an edge connects these two specific vertices, thereby determining whether Al Capone and Albert Anastasia actually share a business relationship.

In a subsequent exercise, students are asked to use PHOEG to hypothesize an upper bound on the \emph{chromatic number} based on the \emph{maximum degree} for \emph{connected graphs}. The pedagogical goal is to guide them toward empirically rediscovering Brooks' theorem~\cite{brooks1941colouring}, or at least closely approaching it. This foundational theorem states that for any connected graph that is neither a \emph{complete graph} nor an \emph{odd cycle}, its chromatic number is at most its maximum degree.

A final exercise asks students to explore general graphs of a given size $m$ that minimize the number of \emph{non-equivalent colorings}---that is, valid vertex colorings that are genuinely distinct and cannot be transformed into one another simply by permuting the color labels. Unlike the previous problems, accurately formulating the conjecture for these presumed extremal graphs is significantly more intricate, primarily because their formal description involves triangular numbers. Students are informed that a conjecture regarding these specific graphs was initially identified in 2014 with the assistance of PHOEG. However, despite continuous efforts to prove or disprove it, it remains an open problem today. Pedagogically, this exercise is extremely valuable. It demonstrates to students the open-ended nature of scientific research and perfectly illustrates a fundamental reality of extremal graph theory: while computational tools like PHOEG can make the empirical identification of a conjecture surprisingly accessible, providing its formal mathematical proof can remain a formidable challenge.

\section{Conclusion}
\label{sec:conclusion}

This paper presented PHOEG, an interactive online tool designed to support both research and education in extremal graph theory. Built around an exact geometrical approach that embeds graphs into a two-dimensional invariant space and computes their convex hulls, PHOEG empowers users to visually explore relationships between graph invariants, identify extremal graphs, and formulate conjectured inequalities.

We detailed the platform's modern architecture, including its comprehensive web interface --- featuring problem definition, interactive polytope visualization, and highly customizable graph rendering --- and its standalone RESTful API for programmatic database access. Furthermore, we demonstrated PHOEG's practical impact in two primary domains. In research, we illustrated its capacity to facilitate the discovery of new conjectures, swiftly identify counterexamples, and provide structural intuition for formal mathematical proofs, effectively serving as a short survey of several results enabled by this geometric approach. In education, we highlighted its successful integration into advanced graph theory courses at the University of Mons and Polytechnique Montr\'eal, where students actively engage with the tool to empirically rediscover classical theorems and grapple with open mathematical problems.

Looking ahead, we plan to continuously expand the underlying database by computing higher graph orders and integrating additional invariants, actively welcoming suggestions from the scientific community. In particular, restricting to specific graph classes such as trees or chemical graphs, for which the combinatorial explosion is significantly more manageable, would allow us to push the database to higher orders.  Ultimately, we hope that PHOEG will continue to evolve as an accessible and invaluable resource for both researchers and educators in the field of graph theory.

\bibliographystyle{acm}
\bibliography{biblio}

\end{document}